\documentclass[12pt]{article}
\usepackage{times}

\newcommand{\doublespace}{
   \renewcommand{\baselinestretch}{1.2}
   \large\normalsize}

\def \Z{\Bbb Z}
\def \C{\Bbb C}

\def \<{\langle}
\def \>{\rangle}

\def \l{\lambda }
\def \L{\Lambda }

\def \pf{\noindent {\bf Proof:} \,}
\def \cg{\chi_g}
\def \cg'{\chi'_g}

\def \d{\delta}

\input amssym.def
\input amssym
\doublespace
\begin{document}
\bibliographystyle{alpha}
\newtheorem{thm}{Theorem}[section]
\newtheorem{thmn}{Theorem}
\newtheorem{prop}[thm]{Proposition}
\newtheorem{cor}[thm]{Corollary}
\newtheorem{lem}[thm]{Lemma}
\newtheorem{rem}[thm]{Remark}
\newtheorem{de}[thm]{Definition}
\newtheorem{hy}[thm]{Hypothesis}
\begin{center}
{\Large {\bf A Spanning Set for VOA Modules}}
\\
\vspace{0.5cm}
Geoffrey Buhl
\\
Department of Mathematics, University of
California, Santa Cruz, CA 95064
\\
{\tt gwbuhl@math.ucsc.edu}
\end{center}

\begin{abstract}
We develop a spanning set for weak modules of $C_2$ co-finite
vertex operator algebras.  This spanning set has finiteness
properties that we use to show weak modules are $C_n$ co-finite
and $A_n(M)$ is finite dimensional.
\end{abstract}

\section{Introduction}

A vertex operator algebra, $V$, is $C_2$ co-finite if the subspace
$\{u_{-2}v: u,v \in V\}$ has finite codimension. This condition,
sometimes referred to as the $C_2$ condition or Zhu's finiteness
condition \cite{Z}, is important in the theory of vertex operator
algebras. Zhu used it, as well as other assumptions, to show
modular invariance of certain trace functions.  One important
feature of the $C_2$ condition is that it is an internal
condition.  Given a vertex operator algebra it is relatively easy
to calculate if $V$ is $C_2$ co-finite. The implications of the
$C_2$ condition are wide ranging, however they are not completely
understood.  The $C_2$ condition implies that $A(V)$ is finite
dimensional, there are a finite number of irreducible admissible
modules \cite{DLM3}, and irreducible admissible modules are
ordinary \cite{KL}. The goal of this paper is to demonstrate new
implications of the $C_2$ condition.  Specifically, we develop new
results about the modules of a vertex operator algebra that
satisfies the $C_2$ condition.  This information sheds new light
on how this internal condition affects the structure of modules.

In recent work, Gaberdiel and Neitzke \cite{GN} develop a spanning
set for vertex operator algebras. They show that $V$ is spanned by
certain vectors of the form $x^1_{-n_1} x^2_{-n_2} \cdots
x^k_{-n_k} {\textbf 1}$, where the modes are strictly decreasing
and less than zero, i.e. $n_1 > n_2 > \cdots
> n_k > 0$. This generating set is finite if $V$ is $C_2$ co-finite.
The principle feature of this spanning set is this no
repeat condition.  Using this result, they prove $C_2$
co-finiteness implies $C_n$ co-finiteness for $n \geq 2$, and the
fusion rules for irreducible admissible modules are finite.

In this paper, we develop an analogous spanning set for modules of
vertex operator algebras.  We will show that under the $C_2$
condition, any weak module generated by $w$ is spanned by certain
elements of the form $x^1_{-n_1} x^2_{-n_2} \cdots x^k_{-n_k}w$.
Here the modes will be decreasing, and each mode will be repeated
at most a finite number of times.   This finite repeat condition,
though not as strong as a no repeat condition, still allows us to
prove a number of results.  With this new spanning set, we extend
the results of Gabediel and Neitzke to modules.  This means we
will demonstrate $A_n(M)$ is finite dimensional and $C_n(M)$ is
finite dimensional, for $M$ finitely generated weak modules.  In a
future paper, we will use this module spanning set to show that
rationality and $C_2$ co-finiteness imply regularity.

The following is a brief preview of the remaining sections of this
paper.  In the second section of this paper, we give the necessary
definitions and notation conventions.  We also present key results
leading up to this paper.  In the third section we develop the
theory of so called singular like vectors. With these vectors, we
are able to reduce expressions with  repeated modes.  In the
fourth section of this paper, we prove our main result which is
the module spanning set with a finite repeat condition. The last
section of this paper gives additional results that quickly follow
from the main theorem.

\section{Preliminaries}

We make the assumption that the reader is somewhat familiar with
the theory of vertex operator algebras (VOAs).  We assume the
definition of a vertex operator algebra as well as some basic
properties. Good reference material is available in papers by
Borcherds \cite{B}; Dong \cite{D}; and Frenkel, Huang, and
Lepowsky \cite{FHL} and in a book by Frenkel, Lepowsky, and Meuman
\cite{FLM}. We begin with some definitions.

\begin{de}
A vertex operator algebra, $V$, is of CFT type if
$V=\bigoplus_{n\geq 0}V(n)$ and $V(0)=span\{ {\textbf 1} \}$.
\end{de}

Throughout out this paper, we assume $V$ is of CFT type.

\begin{de}
For $V$ a VOA, $C_n(V)=\{v_{-n}w \mid v,w \in V\}$
\end{de}

\begin{de}
$V$, a VOA,  is called $C_n$ co-finite if $V \slash C_n(V)$ is
finite dimensional.  For $n=2$,this is often called the $C_2$
condition.
\end{de}

$C_2$ co-finiteness is an important assumption in Zhu's work
demonstrating modularity of certain functions.  This paper is
about module spanning sets, and since there are a few different
flavors of VOA modules, we now define weak, admissible, and
ordinary modules.

\begin{de}
A weak $V$ module is a vector space $M$ with a linear map \\
$Y_M:V \rightarrow End(M)[[z,z^{-1}]]$ where $v \mapsto
Y_M(v,z)=\sum_{n \in \Z}v_n z^{-n-1}$, $v_n \in End(M)$. In
addition $Y_M$ satisfies the following:

1) $v_nw=0$ for $n>>0$ where $v \in V$ and $w \in M$

2) $Y_M( {\textbf 1},z)=Id_M$

3) The Jacobi Identity:
\begin{eqnarray}
z_0^{-1}\d ({z_1 - z_2 \over z_0})Y_M(u,z_1)Y_M(v,z_2)-
z_0^{-1} \d ({z_2- z_1 \over -z_0})Y_M(v,z_2)Y_M(u,z_1) \nonumber \\
=z_2^{-1} \d ({z_1- z_0 \over z_2})Y_M(Y(u,z_0)v,z_2)
\end{eqnarray}
\end{de}

There are two important consequences of this definition.  Weak
modules admit a Virasoro representation under the action of
$\omega$, the Virasoro vector.  Also weak modules satisfy the
$L(-1)$ derivation property. The distinctive feature of weak
modules is that they have no grading. Admissible and ordinary
module are both graded.

\begin{de}
An admissible $V$ module is a weak $V$ module which carries a
nonnegative integer grading, $M=\bigoplus_{n \geq 0} M(n)$, such
that if $v \in V(r)$ then $v_m M(n) \subseteq M(n+r-m-1)$
\end{de}

So for an admissible module, we have added a grading with a bottom
the level, and the action of $V$ respects the grading.

\begin{de}
An ordinary $V$ module is a weak $V$ module which carries a $\C$
grading, $M=\bigoplus_{\l \in \C} M_{\l}$, such that:

1) $dim(M_{\l})< \infty$

2) $M_{\l+n=0}$ for fixed $\l$ and $n<<0$

3) $L(0)w=\l w=wt(w) w$, for $w \in M$
\end{de}

Although the definition of a $\C$ grading may seem weaker than a
$\Z$ grading, the requirement that each graded piece of an
ordinary modules must be finite dimensional is a strong condition.
It turns out that any ordinary module is admissible. So we have
this set of inclusions:

$$\{\mbox{ordinary modules}\} \subseteq \{\mbox{admissible modules}\} \subseteq \{\mbox{weak modules}\}$$

In addition to the above definitions, we need to refer to results
by Gaberdiel and Neitzke \cite{GN}, and their work in determining
a spanning set for vertex operator algebras.  There are three
pertinent results of theirs that are explained below.  First we
describe the generating set. Let $\{\bar{x}^i\}_{i \in I}$ be a
basis of $V \slash C_2(V)$, where $\bar{x}^i=x^i + C_2(V)$, and
$x^i$ is a homogenous vector. So $\bar{X}=\{x^i\}_{i \in I}$ is a
set of elements in $V$ which are representatives of a basis for $V
\slash C_2(V)$.

\begin{thm}\cite{GN}
Let $V$ be a vertex operator algebra, then $V$ is spanned by
elements of the form
$$x^{1}_{-n_1} x^{2}_{-n_2} \cdots x^{k}_{-n_k} {\textbf 1}$$
where $n_1>n_2> \cdots >n_k > 0$ and $x^{j} \in \bar{X}$ for $1
\leq j \leq k$.
\end{thm}

Throughout this paper will refer to the elements described in this
theorem as VOA spanning set elements.  This is in contrast to the
module spanning set elements that are the goal of this paper. This
theorem tells us the modes in this spanning set are strictly
decreasing, and this is what we mean by a no repeat condition.
This spanning set is especially usefully when we look at the $C_n$
spaces.  In fact the next result follows very quickly from this no
repeat condition.

\begin{thm} \cite{GN}
Suppose $V \slash C_2(V)$ is finite dimensional, then $V \slash
C_n(V)$ is finite dimensional for $n \geq 2$.
\end{thm}

In their paper, Gaberdiel and Neitzke also formulate a spanning
set for modules.  Unfortunately, this formulation does not have a
no repeat condition.  Their repetition restriction is in terms of
the weights of the modes. In their module spanning set, the weight
of the modes are decreasing, but the inequality is not strict,
which means that a particular mode could be repeated indefinitely.
The following is the definition of the weight of the mode.

\begin{de}
For $u \in V$, a homogeneous vector, and $n \in \Z$,
$wt(u_n)=wt(u)-n-1$.
\end{de}

\begin{thm} \cite{GN}
Let $M$ be an admissible  $V$ module.  Then $M$ is spanned by
elements of the form
$$x^1_{-n_1} x^2_{-n_2} \cdots x^k_{-n_k} u$$
where u is a lowest weight vector and $wt(x^1_{-n_1}) \leq
wt(x^2_{-n_2}) \leq \cdots \leq wt(x^k_{-n_k}) \leq 0$ and $x^{j}
\in \bar{X}$ for $1 \leq j \leq k$.
\end{thm}

The goal of this paper is a module spanning set analogous to the
VOA spanning set result by Gaberdiel and Neitzke under the
additional condition of $C_2$ co-finiteness.  In this new spanning
set, modes are only allow to repeat an finite number of times,
which is slightly weaker that the condition that modes can only be
repeated once. It turns out however that this finite repeat
condition is sufficient to demonstrate two nice finiteness
properties for weak modules: $A_n(M)$ finite dimensionality and
$C_n(M)$ co-finiteness.

In this paper the term, mode, is abused slightly.  The term, mode,
usually refers to an endomorphism, $u_n$, where $u \in V$ and $n
\in \Z$. Sometimes, in this paper it refers to the indexing number
$n$. For example, in the term, nonnegative mode,``mode'' refers to
a mode of the form $u_n$ with $n \geq 0$. It should be apparent
what ``mode'' is referring to by the context.

\section{Singular Like Vectors}

    In order to limit the number of repeated modes, we need a method
for shortening spanning set elements that have too many
repetitions of a certain mode.  To meet this end, we must develop
the theory of so called singular like vectors.  The reason for
this name is because the vectors described in this section are
reminiscent of singular vectors in the Virasoro algebra \cite{FF}.

The goal of this section is twofold.  First, we choose vectors of
the form $x^{1}_{-1} x^{2}_{-1} \cdots x^{k}_{-1}{\textbf 1}$ in our
VOA and rewrite them as a sum of spanning set elements,
$x^{r_1}_{-n_{r_1}} x^{r_2}_{-n_{r_2}} \cdots
x^{r_l}_{-n_{r_l}}{\textbf 1}$ where $l<k$. Second, we calculate the
vertex operators of these singular like vectors. Isolating certain
coefficients of these vertex operators is key to limiting the
repetition of modes in the module spanning set that is the main
result of this paper.

Henceforth, we are working under the assumption that $V$, our VOA,
is $C_2$ co-finite.  So now, $\bar{X}$ is a finite set homogenous
of elements in $V$ which are representatives of a basis for $V
\slash C_2(V)$. We know that $V$ is spanned by elements
of the form $x^{1}_{-n_1} x^{2}_{-n_k} \cdots x^{k}_{-n_k} {\textbf
1}$ where $x^{i} \in \bar{X}$ for $1 \leq i \leq k$ and the modes
are strictly decreasing.

We can simplify this $\bar{X}$ slightly.  The vacuum, ${\textbf 1}$,
is not an element of $C_2(V)$ so we could choose a basis $\bar{X}$
such that ${\textbf 1} \in \bar{X}$.  The only mode for ${\textbf 1}$
which is nonzero is ${\textbf 1}_{-1}$ but this is the identity
endomorphism. If we define $X= \bar{X}-\{{\textbf 1}\}$, the results
of Gaberdiel and Neitzke still hold, but $X$ is one element
smaller.  Note that this means the minimum weight of any vector in
$X$ is $1$.

The first set of results in this section establish that we can
rewrite a certain type of vector as a sum of VOA spanning elements
of strictly shorter length.  The proof involves comparing the
weights of vectors.  So we start the following definition.

\begin{de}
Let $B=\max_{x \in X} \{ wt(x) \}$.
\end{de}

This means $B$ is the largest weight of any vector not in $C_2(V)$

\begin{lem}

Let $x^{i} \in X$ for $1 \leq i \leq k$, $k$ a positive integer.

\begin{equation}
wt(x^{1}_{-1} x^{2}_{-1} \cdots x^{k}_{-1}{\textbf 1}) \leq Bk
\label{minusone}
\end{equation}

\end{lem}

\pf
\begin{eqnarray}
wt(x^{1}_{-1} x^{2}_{-1} \cdots x^{k}_{-1}{\textbf 1}) &=&
\sum^k_{i=1} wt(x^{i}) \\
&\leq& k \max_{x \in X} \{wt(x)\} \\
&=& Bk
\end{eqnarray}
\mbox{$\square$}

Now that we have a maximum weight for a vector of the form
$x^1_{-1} x^2_{-1}\cdots x^n_{-1}{\textbf 1}$, we compare its weight
to a spanning set vectors weight. The next lemma will tell us the
minimum weight of a VOA spanning set vector
 of length $l$.

\begin{lem}
Let $x^{i} \in X$ for $1 \leq i \leq l$, $l$ a positive integer.
If $n_1 > n_2 > \cdots > n_l > 0$ then,
\begin{equation}
wt(x^{1}_{-n_1} x^{2}_{-n_2} \cdots x^{l}_{-n_l}{\textbf 1}) \geq
{l(l+1) \over 2}
\end{equation}
\end{lem}

\pf
\begin{eqnarray}
wt(x^{1}_{-n_1} x^{2}_{-n_2} \cdots x^{l}_{-n_l}{\textbf 1}) &=&
\sum^l_{i=1} wt(x^{i}) + n_i - 1 \\
&\geq& \sum^l_{t=1} n_i \\
&\geq& \sum^l_{t=1} i \\
&=& {l(l+1) \over 2}
\end{eqnarray}
$\square$

Note here that the minimum weight of a VOA spanning set element
increases quadratically as a function of length, while the maximum
weight of a vector of the form $$x^{1}_{-1} x^{2}_{-1} \cdots
x^{k}_{-1} {\textbf 1}$$ increases linearly as a function of length.
Using the previous two lemmas, we will now show that if we have a
vector with a sufficient number of $-1$ modes, it can be rewritten
as a sum of spanning set elements of strictly shorter length.

\begin{lem}
\label{lemma33} Let $x^{i} \in X$ for $1 \leq i \leq k$ and $y^{j}
\in X$ for $1 \leq j \leq l$. If $2B-1 < k \leq l$ and $n_1 > n_2
> \cdots > n_l > 0$, then
\begin{equation}
wt(x^{1}_{-1} x^{2}_{-1} \cdots x^{k}_{-1}{\textbf 1}) <
wt(y^{1}_{-n_1} y^{2}_{-n_2} \cdots y^{l}_{-n_l}{\textbf 1})
\end{equation}
\end{lem}

\pf
From the previous lemmas we have,
$$wt(x^{1}_{-1} x^{2}_{-1} \cdots x^{k}_{-1}{\textbf 1}) \leq Bk$$
and
$$wt(y^{1}_{-n_1} y^{2}_{-n_2} \cdots y^{l}_{-n_l}{\textbf 1}) \geq
{l(l+1) \over 2}$$
If  $l \geq k$ then,
\begin{equation}
{l(l+1) \over 2} \geq {k(k+1) \over 2}
\end{equation}
Since $k>2B-1$, then ${k(k+1) \over 2} > Bk$.  Thus,
$wt(y^{1}_{-n_1} y^{2}_{-n_2} \cdots y^{l}_{-n_l}{\textbf 1})>Bk$ ,
and we have, finally, that $wt(x^{1}_{-1} x^{2}_{-1} \cdots
x^{k}_{-1}{\textbf 1}) < wt(y^{1}_{-n_1} y^{2}_{-n_2} \cdots
y^{l}_{-n_l}{\textbf 1})$.
\mbox{$\square$}

Before we continue, we need to recall a basic fact about $C_2$
co-finite vertex operator algebras.

\begin{rem}
If $V$ is $C_2$ co-finite, and $V=\bigoplus_{i \geq 0}V_i$ is the
weight space decomposition of $V$.  Then for some $N>0$,
$\bigoplus_{i\geq N} V_i \subset C_2(V)$.
\end{rem}

\begin{de}
Let $Q=\max\{N,2B-1\}+1$.
\end{de}

This $Q$ will play an important role in the proof of the module
spanning set.  We start out, however, by showing that certain
vectors long enough can be rewritten in terms of shorter vectors.
In particular, if we have a vector composed of the product of $Q$,
minus one modes acting on the vacuum, we can rewrite it in terms
of VOA spanning set elements which only have at most $Q-1$ modes.

\begin{prop}
\label{37}
If $k \geq Q$,
$$x^{1}_{-1} x^{2}_{-1} \cdots x^{k}_{-1}{\textbf 1}= \sum_{r
\in R} x^{r_1}_{-n_{r_1}} x^{r_2}_{-n_{r_2}} \cdots
x^{r_l}_{-n_{r_l}}{\textbf 1}$$ with $l<k$ and where $x^{i} \in X$ for
$1 \leq i \leq k$; $n_{r_1} > n_{r_2} > \cdots > n_{r_l} > 0$ for
fixed $r$; $x^{r_t} \in X$ for $1 \leq t \leq l$; and $R$ a finite
index set.
\end{prop}

\pf If $k>N$, then $x^{1}_{-1} x^{2}_{-1} \cdots x^{k}_{-1}{\textbf 1}
\in C_2(V)$,  and we can write

\begin{equation}
\label{slv} x^{1}_{-1} x^{2}_{-1} \cdots x^{k}_{-1}{\textbf 1}=
\sum_{r \in R} x^{r_1}_{-n_{r_1}} x^{r_2}_{-n_{r_2}} \cdots
x^{r_l}_{-n_{r_l}}{\textbf 1}
\end{equation}
where $n_1 \geq 2$ and $n_1 > n_2 > \cdots > n_l > 0$.  So we have
rewritten $x^{1}_{-1} x^{2}_{-1} \cdots x^{k}_{-1}{\textbf 1}$ as sum
of spanning set elements in $C_2(V)$.  Now assume to the contrary,
that for all $k$ and $r$, $l \geq k$. Now by the lemma
\ref{lemma33}, if $k>2B-1$,

$$wt(x^{1}_{-1} x^{2}_{-1} \cdots x^{k}_{-1}{\textbf 1}) <
wt(x^{r_1}_{-n_{r_1}} x^{r_2}_{-n_{r_2}} \cdots
x^{r_l}_{-n_{r_l}}{\textbf 1})$$
This is a contradiction since
$$wt(x^{1}_{-1} x^{2}_{-1} \cdots x^{k}_{-1}{\textbf 1}) =
wt(x^{r_1}_{-n_{r_1}} x^{r_2}_{-n_{r_2}} \cdots
x^{j_l}_{-n_{r_l}}{\textbf 1})$$
 Thus, if $k \geq Q$,
$l<k$ for all $r$.

\mbox{$\square$}

We call vectors of the form (\ref{slv}) singular like because they
are similar to singular vectors in the Virasoro algebra. The main
point of this lemma is that there is a uniform length for which
any product of $-1$ modes of that length can be rewritten as sum
of strictly shorter spanning set elements.

The next step is to calculate the vertex operators of these
singular like vectors.  This will allow us later on to derive an
identity that enables us to shorten repeated modes.  In order to
calculate the vertex operators we need to recall a formula.

\begin{rem}
\begin{equation}
Y(u_n v,z)=Res_{z_1} \{ (z_1-z)^n Y(u,z_1) Y(v,z) - (-z+z_1)^n Y(v,z) Y(u,z_1) \}
\end{equation}

By taking the residue, we obtain

\begin{equation}
Y(u_n v,z)= \sum\limits_{m \geq 0} (-1^m) {{n} \atopwithdelims ()
m} z^m u_{n-m} Y(v,z) - \sum\limits_{m \geq 0} (-1)^{n-m} {{n}
\atopwithdelims () m} z^{n-m}Y(v,z) u_{m}
\end{equation}

In particular, if $n=-1$, we have

\begin{equation}
Y(u_{-1} v,z)=\sum\limits_{m \geq 0} z^m u_{-1-m} Y(v,z) +
\sum\limits_{m \geq 0} z^{-1-m}Y(v,z) u_{m} \label{calcone}
\end{equation}

\end{rem}

To perform the next calculation, we will need to apply the above
formula Q times to the left hand side of our singular like vector.
But before we can do the next calculation, we need to describe a
certain indexing set.  This set, in its ordering, will describe
how operators are rearranged when we apply (\ref{calcone})
multiple times.

\begin{de}
$\L^i_n=\{ (k_1, k_2, \ldots, k_i) : \{k_1, k_2, \ldots, k_i\}
\subseteq \{1, 2, \ldots, n\}, k_1 < k_2 < \cdots <k_i\}$.  Now if
$\l \in \L^i_n$ then $\l_j$ is the $j^{th}$ element of $\l$.
\end{de}

So what we have is an $i$ element subset of $\{1, 2, \ldots, n\}$
where the order is specified. The elements are placed increasing
order.  $\L^i_n$ is then the set of these $i$ element subsets.

\begin{de}
If $\l \in \L^i_n$ then $\bar{\l}=(k_{i+1}, k_{i+2}, \ldots,
k_{n})$ where $\{k_{i+1}, k_{i+2}, \ldots, k_{n}\}$ is the
compliment of $\{k_{1}, k_{2}, \ldots, k_{i}\}$ in $\{1, 2,
\ldots, n\}$ and  $k_{i+1} > k_{i+2} > \cdots > k_{n}$.  Also
$\bar{\l}_j$ is the $j^{th}$ element of $\bar{\l}$.
\end{de}

So in the compliment to $\l$, the order is reversed.  The elements
are in decreasing order.   For example, if $n=8$, $\L^i_8$ is the
set of $i$ element subsets of $\{1,2,3,4,5,6,7,8\}$.  If $i=3$,
let $\l=\{2,5,8\}$, then $\bar\l=\{7,6,4,3,1\}$.  If $i=4$ and
$\l=\{3,4,6,7\}$, then $\bar\l=\{8,5,2,1\}$.

\begin{prop}
\label{vertexopofsv} Let $x^1, \ldots, x^n, v \in V$, then
\begin{eqnarray}
\lefteqn{Y(x^{1}_{-1} x^{2}_{-1} \cdots x^{n}_{-1}v,z) =} \nonumber \\
&&\sum_{i=0}^{n} \sum_{\l \in \L^i_n} \sum_{m_{i} \geq 0}
(\prod_{j=1}^{i} x^{\l_j}_{-1-m_{\l_j}}z^{m_{\l_j}}) Y(v,z)
(\prod_{j=i+1}^{n}
x^{\bar{\l}_j}_{m_{\bar\l_j}}z^{-1-m_{\bar\l_j}})
\end{eqnarray}

\end{prop}

\pf

We proceed by induction.  For the case where $n=1$, we use
(\ref{calcone}). When $n=1$ we have $\L^0_1=\{ \emptyset \}$ and
$\L^1_1=\{ \{1\} \}$.  So we have
\begin{eqnarray}
Y(x_{-1} v,z)&=&\sum\limits_{m \geq 0} z^m x_{-1-m} Y(v,z) +
\sum\limits_{m \geq 0} z^{-1-m} Y(v,z) x_{m}\\
&=& \sum_{i=0}^1 \sum_{m \geq 0} (\prod_{j=1}^i
x_{-1-m}z^m)(Y(v,z)) (\prod_{j=i+1}^1 x_m z^{-1-m})
\end{eqnarray}

Now assume true for $n-1$, then

\begin{eqnarray}
\lefteqn{Y(x^{1}_{-1} x^{2}_{-1} \cdots x^{n}_{-1}v,z)} \\
&=& \sum_{i=0}^{n-1} \sum_{\l \in \L^i_{n-1}} \sum_{m_{i} \geq
0} (\prod_{j=1}^{i} x^{\l_j}_{-1-m_{\l_j}} z^{m_{\l_j}})
(Y(x^n_{-1}v,z))
(\prod_{j=i+1}^{n-1} x^{\bar{\l}_j}_{m_{\l_j}}z^{-1-m_{\bar\l_j}}) \\
&=& \sum_{i=0}^{n-1} \sum_{\l \in \L^i_{n-1}} \sum_{m_{i} \geq
0}
(\prod_{j=1}^{i} x^{\l_j}_{-1-m_{\l_j}} z^{m_{\l_j}}) \\
&& \cdot (\sum\limits_{m_n \geq 0} z^{m_n} x_{-1-m_n} Y(v,z) +
\sum\limits_{m_n \geq 0} z^{-1-m_n} Y(v,z) x_{m_n}) \\
&& \cdot (\prod_{j=i+1}^{n-1} x^{\bar{\l}_j}_{m_{\bar\l_j}}z^{-1-m_{\bar\l_j}}) \\
&=&\sum_{i=0}^{n} \sum_{\l \in \L^i_n} \sum_{m_{i} \geq 0}
(\prod_{j=1}^{i} x^{\l_j}_{-1-m_{\l_j}}z^{m_{\l_j}}) Y(v,z)
(\prod_{j=i+1}^{n}
x^{\bar{\l}_j}_{m_{\bar\l_j}}z^{-1-m_{\bar\l_j}})
\end{eqnarray}

\mbox{$\square$}

So, in particular when $v={\textbf 1}$, we have

\begin{eqnarray}
\label{svvo}
\lefteqn{Y(x^{1}_{-1} x^{2}_{-1} \cdots x^{n}_{-1} {\textbf 1},z) =} \nonumber \\
&&\sum_{i=0}^{n} \sum_{\l \in \L^i_n} \sum_{m_i \geq 0}
(\prod_{j=1}^{i}
x^{\l_j}_{-1-m_{\l_j}}z^{m_{\l_j}})(\prod_{j=i+1}^{n}
x^{\bar{\l}_j}_{m_{\bar\l_j}}z^{-1-m_{\bar\l_j}})
\end{eqnarray}

Now we consider $x^1, \ldots, x^Q \in X$.  By Proposition
\ref{37},

$$x^{1}_{-1} x^{2}_{-1} \cdots x^{Q}_{-1}{\textbf 1}=
\sum_{r \in R} x^{r_1}_{-n_{r_1}} x^{r_2}_{-n_{r_2}} \cdots
x^{r_l}_{-n_{r_l}}{\textbf 1}$$

\noindent where  $l<Q$; $n_{r_1}
> n_{r_2} > \cdots > n_{r_l} > 0$ for fixed $r$; $x^{r_t} \in X$ for $1 \geq t \geq
l$; and $R$ a finite index set.  Substituting $\sum_{r \in R}
x^{r_1}_{-n_{r_1}} x^{r_2}_{-n_{r_2}} \cdots
x^{r_l}_{-n_{r_l}}{\textbf 1}$ for $x^{1}_{-1} x^{2}_{-1} \cdots
x^{Q}_{-1}{\textbf 1}$ in (\ref{svvo}), we get:

\begin{eqnarray}
\label{bloom} \lefteqn{Y(\sum_{r \in R}
x^{r_1}_{-n_{r_1}} x^{r_2}_{-n_{r_2}} \cdots x^{r_l}_{-n_{r_l}}{\textbf 1},z) =} \nonumber \\
&&\sum_{i=0}^{Q} \sum_{\l \in \L^i_Q} \sum_{m_{i} \geq 0}
(\prod_{j=1}^{i}
x^{\l_j}_{-1-m_{\l_j}}z^{m_{\l_j}})(\prod_{j=i+1}^{Q}
x^{\bar{\l}_j}_{m_{\bar\l_j}}z^{-1-m_{\bar\l_j}})
\end{eqnarray}

\noindent where $R$ is a finite index set, $x^1, \cdots, x^Q,
x^{r_1}, \ldots, x^{r_l} \in X$, $n_{r_1}
> n_{r_2} > \cdots > n_{r_l} > 0$ for fixed $r$, and $l <Q$.

\section{A Module Spanning Set}

To reiterate our setting $V=(V,Y,{\textbf 1},\omega)$ is a $C_2$
co-finite vertex operator algebra of CFT type.    In this section
we prove the main theorem of this paper.  This theorem will give a
spanning set for weak $V$ modules using the set ${X}$. This
spanning set will have mode repetition restrictions similar to
those given by Gaberdiel and Neitzke for vertex operator algebras.
We begin by formulating a few lemmas, which give us identities
that we use to impose mode repetition restrictions on the module
spanning set elements.

\begin{rem}  Borcherds's Identity \cite{B}: Let $u,v \in V$ and $k,q,r \in \Z$
\begin{eqnarray}
\lefteqn{\sum_{i \geq 0} {{-k} \atopwithdelims () i} (u_{-r+i}v)_{-k-q-i}}\\
&=&\sum_{i \geq 0} (-1)^i {{-r} \atopwithdelims () i}
\{u_{-k-r-i}v_{-q+i} - (-1)^{-r} v_{-q-r-i} u_{-k+i}\}
\end{eqnarray}
This is the component form of the Jacobi identity.
\end{rem}

We first recall two formulas resulting from Borcherds's Identity.

\begin{lem}
\label{41} Let $u,v \in V$ and $k,q \in \Z$, then
$$[ u_{-k}, v_{-q}]=\sum\limits_{i \geq 0}
{{-k} \atopwithdelims () i} (u_{i}v)_{-k-q-i}$$
\end{lem}

\pf  In Borcherds's Identity, let $r=0$.\mbox{$\square$}

\begin{lem}
\label{42}
Let $u,v \in V$ and $r,q \in \Z$, then
\begin{eqnarray*}
(u_{-r}v)_{-q}&=&\sum\limits_{i \geq 0}(-1)^i
{{-r} \atopwithdelims () i}u_{-r-i}v_{-q+i}\\
&&-\sum\limits_{i \geq 0}(-1)^{i-r} {{-r} \atopwithdelims ()
i}v_{-r-q-i}u_{i}
\end{eqnarray*}
\end{lem}

\pf In Borcherds's Identity, let $k=0$. \mbox{$\square$}

It is worth noting that in the previous two lemmas the weight of
operators on either sides of the equations are equal.  This can be
verified by simple calculations.  The third formula is a special
case of the previous lemma.

\begin{lem}
\label{43} Let $u,v \in V$ and $n \in \Z$, then
\begin{eqnarray*}
u_{-n} v_{-n} = (u_{-1}v)_{-2n+1}- \sum\limits_{i< 0}u_{-1-i}v_{-2n+1+i}
+ \sum\limits_{i \geq 0}v_{-2n-i} u_{i}
\end{eqnarray*}
\end{lem}

\pf This follows from the previous lemma with $q=2n-l$ and $r=1$.
\mbox{$\square$}

This third lemma will be the identity we use to limit repeated
negative modes.

\begin{de}
Let $W$ be a weak $V$ module, and $w \in W$. For each $x \in X$
there exists $l_x \geq 0$ such that $x_{l_x}w \ne 0$ but
$x_{l_x+k}w=0$ for all $k>0$.  Now define $L= \max_{x \in X}
\{l_x\}+1$.
\end{de}

When we are examining expressions involving products of modes, we
do not need to look at modes $L$ or larger, because they will kill
the given $w$.

In addition to Lemma \ref{43}, we need another lemma that allow us
to impose repetition restrictions on nonnegative modes. To obtain
this new formula, we take a reside of equation in Proposition
\ref{vertexopofsv}. The goal is to isolate certain coefficients of
$z$ where we have repeated positive modes. We will eventually find
an expression for endomorphisms of the form $x^{1}_n \cdots
x^{Q}_n$ where $0 \leq n \leq L$.

\begin{lem}
Let $x^{1}_{-1}\cdots x^{Q}_{-1}{\textbf 1}= \sum_{r \in R}
x^{r_1}_{-n_{r_1}} x^{r_2}_{-n_{r_2}} \cdots
x^{r_l}_{-n_{r_l}}{\textbf 1}$ where $x^{i}, x^{r_t} \in X$ for $1
\leq i \leq Q$ and $1 \leq t \leq l$; $l<Q$; and $n_{r_1}
> n_{r_2} > \cdots > n_{r_l} > 0$ for fixed r.  Then for $0 \leq k
\leq L$,

\begin{eqnarray}
\lefteqn{z^{Q(-1+k-L)}x^Q_{L-k} x^{Q-1}_{L-k} \cdots x^1_{L-k}} \\
&=&Y(\sum_{r \in R}
x^{r_1}_{-n_{r_1}} x^{r_2}_{-n_{r_2}} \cdots x^{r_l}_{-n_{r_l}}{\textbf 1},z) \\
&&- \sum_{i=1}^{Q} \sum_{\l \in \L^i_Q} \sum_{m_{i} \geq 0}
(\prod_{j=1}^{i}
x^{\l_j}_{-1-m_{\l_j}}z^{m_{\l_j}})(\prod_{j=i+1}^{Q}
x^{\bar{\l}_j}_{m_{\bar\l_j}}z^{-1-m_{\bar\l_j}}) \\
&&- \sum_{m_{j} \geq 0, 1 \leq j \leq Q} z^{-Q-\sum^Q_{j=1} m_{j}}
(x^{Q}_{m_{Q}}\cdots x^{1}_{m_{1}}) \label{46s3}
\end{eqnarray}
where, in (\ref{46s3}), for at least one $j$, $m_{j} \neq L-k$.
\end{lem}

\pf From (\ref{bloom}), we have

\begin{eqnarray}
\lefteqn{Y(\sum_{r \in R} x^{r_1}_{-n_{r_1}} x^{r_2}_{-n_{r_2}}
\cdots x^{r_l}_{-n_{r_l}}{\textbf 1},z)=} \nonumber \\
&&\sum_{i=0}^{Q} \sum_{\l \in \L^i_Q} \sum_{m_{i} \geq 0}
(\prod_{j=1}^{i}
-x^{\l_j}_{-1-m_{\l_j}}z^{m_{\l_j}})(\prod_{j=i+1}^{Q}
(x^{\bar{\l}_j}_{m_{\bar\l_j}}z^{-1-m_{\bar\l_j}}) \label{vosv}
\end{eqnarray}

From this we will isolate the mode of the form $x^Q_{L-k} \cdots
x^1_{L-k}$ where $L$ is defined as above, and $0 \geq k > L$. That
is, $x_{L-k}$ is a nonnegative mode.

First observe that all nonnegative modes occur when $i=0$.  If $i
\neq 0$ we will have at least one mode of the form
$x^{\l_j}_{-1-m_{\l_j}}$. Now we examine the summand on the right
hand side of (\ref{vosv}) when $i=0$.  Note that when $i=0$, there
is only one possible $\l$, and it is the empty set.  So
$\bar{\l}=\{Q,Q-1, \ldots, 1\}$.  An expression with only positive
modes has the form:

$$\sum_{m_{j} \geq 0, 1 \leq j \leq Q} z^{-Q-\sum^Q_{j=1}
m_{\l_j}} x^{Q}_{m_Q}x^{Q-1}_{m_{Q-1}} \cdots x^{1}_{m_1}$$ Now we
are able to isolate  $x^Q_{L-k} \cdots x^1_{L-k}z^{Q(-1+k-L)}$,
using (\ref{vosv}).

\begin{eqnarray}
\lefteqn{z^{Q(-1+k-L)}x^Q_{L-k} x^{Q-1}_{L-k} \cdots x^1_{L-k}} \\
&=&Y(\sum_{r \in R}
x^{r_1}_{-n_{r_1}} x^{r_2}_{-n_{r_2}} \cdots x^{r_l}_{-n_{r_l}}{\textbf 1},z) \\
&&- \sum_{i=1}^{Q} \sum_{\l \in \L^i_Q} \sum_{m_{i} \geq 0}
(\prod_{j=1}^{i}
x^{\l_j}_{-1-m_{\l_j}}z^{m_{\l_j}})(\prod_{j=i+1}^{Q}
x^{\bar{\l}_j}_{m_{\bar\l_j}}z^{-1-m_{\bar\l_j}}) \\
&&- \sum_{m_{j} \geq 0, 1 \leq j \leq Q} z^{-Q-\sum^Q_{j=1} m_{j}}
(x^{Q}_{m_{Q}}\cdots x^{1}_{m_{1}}) \label{p46s2}
\end{eqnarray}
where in (\ref{p46s2}) for at least one $j$, $m_{j} \neq L-k$.
This means we have isolated the only term for which all modes are
$L-k$. In (\ref{47s2}), there may be several modes equal to $L-k$,
but at least one of the modes is not $L-k$. \mbox{$\square$}

The next step is to take a residue of this result to isolate
$x_{L-k}^{Q} x_{L-k}^{{Q-1}} \cdots x_{L-k}^{1}$.  This will be
our new identity that we will use to impose repetition
restrictions on nonnegative modes

\begin{lem}
\label{mainlemma} Let $x^{1}_{-1}\cdots x^{Q}_{-1}{\textbf 1}= \sum_{r
\in R} x^{r_1}_{-n_{r_1}} x^{r_2}_{-n_{r_2}} \cdots
x^{r_l}_{-n_{r_l}}{\textbf 1}$ where $x^{i}, x^{r_t} \in X$ for $1
\leq i \leq Q$ and $1 \leq t \leq l$, $l<Q$, and $n_{r_1}
> n_{r_2} > \cdots > n_{r_l} > 0$ for fixed r.  Then

\begin{eqnarray}
\lefteqn{x_{L-k}^{Q} x_{L-k}^{{Q-1}} \cdots x_{L-k}^{1} } \\
&=&Res_z \{ Y(\sum_{r \in R} x^{r_1}_{-n_{r_1}} x^{r_2}_{-n_{r_2}}
\cdots x^{r_l}_{-n_{r_l}}{\textbf 1},z)
 z^{Q(-L-1+k)-1}\}  \\
&&- \sum_{i=1}^{Q} \sum_{\l \in \L^i_Q} \sum_{m_{i} \geq 0}
(\prod_{j=1}^{i} x^{\l_j}_{-1-m_{\l_j}})(\prod_{j=i+1}^{Q}
x^{\bar{\l}_j}_{m_{\bar\l_j}}) \label{47s2}\\
&&- \sum_{m_{j} \geq 0, 1 \leq j \leq Q} x^{Q}_{m_{Q}}
x^{Q-1}_{m_{Q-1}}\cdots x^{1}_{m_{1}} \label{47s3}
\end{eqnarray}

Where in (\ref{47s2}), $\sum_{j=1}^{i}(-1-m_{\l_j}) +
\sum_{j=i+1}^{Q}(m_{\bar\l_j}) = Q(L-k)$, and in (\ref{47s3}),
$\sum_{j=1}^{Q}m_{j}=Q(L-k)$ and $m_{j} \neq L-k$ for some $j$.
\end{lem}

When we take the residue, we impose some restrictions on the modes
on the right hand side of this identity.  These restrictions are
that the sum of the modes on the left hand side of the equation
equals the sum of the modes on the left hand side of the equation.
Since the the vectors $x$ are the same on both sides, this means
the the weight of the operators on both sides of the equations are
equal.

\pf We multiply the equation in the statement of Lemma \ref{42} by
$z^{Q(1-k+L)-1}$ and take the residue to obtain the following:

\begin{eqnarray}
\lefteqn{Res_z\{z^{-1}x^Q_{L-k} x^{Q-1}_{L-k} \cdots x^1_{L-k}\}} \\
&=& Res_z\{z^{Q(1-k+L)-1}Y(\sum_{r \in R}
x^{r_1}_{-n_{r_1}} x^{r_2}_{-n_{r_2}} \cdots x^{r_l}_{-n_{r_l}}{\textbf 1},z) \\
&&- z^{Q(1-k+L)-1} \sum_{i=1}^{Q} \sum_{\l \in \L^i_Q} \sum_{m_{i}
\geq 0} (\prod_{j=1}^{i} x^{\l_j}_{-1-m_{\l_j}} z^{m_{\l_j}})
(\prod_{j=i+1}^{Q} x^{\bar{\l}_j}_{m_{\bar\l_j}} z^{-1-m_{\bar\l_j}}) \\
&&- z^{Q(1-k+L)-1} \sum_{m_{j} \geq 0, 1 \leq j \leq Q}
z^{-Q-\sum^Q_{j=1} m_{j}} (x^{Q}_{m_{Q}} \cdots x^{1}_{m_{1}})\}
\label{p147s2}
\end{eqnarray}
Where in (\ref{p147s2}), for at least one $j$, $m_{j} \neq L-k$.
Evaluating the residue we have,

\begin{eqnarray}
\lefteqn{x^Q_{L-k} x^{Q-1}_{L-k} \cdots x^1_{L-k}} \\
&=& Res_z\{z^{Q(1-k+L)-1}Y(\sum_{r \in R}
x^{r_1}_{-n_{r_1}} x^{r_2}_{-n_{r_2}} \cdots x^{r_l}_{-n_{r_l}}{\textbf 1},z)\} \\
&&- \sum_{i=1}^{Q} \sum_{\l \in \L^i_Q} \sum_{m_{i} \geq 0}
(\prod_{j=1}^{i} x^{\l_j}_{-1-m_{\l_j}})(\prod_{j=i+1}^{Q}
x^{\bar{\l}_j}_{m_{\bar\l_j}}) \label{p247s1}\\
&&- \sum_{m_{j} \geq 0, 1 \leq j \leq Q} (x^{Q}_{m_{Q}}\cdots
x^{1}_{m_{1}}) \label{p247s2}
\end{eqnarray}

Where in (\ref{p247s2}),for at least one $j$, $m_{j} \neq L-k$. By
taking the residue we impose the restriction
$Q(1-k+L)-1-Q-\sum_{j=1}^{Q}m_{j} =-1$.  So, we have the
condition:

$$\sum_{j=1}^{Q}m_{j}=Q(L-k)$$

In addition in (\ref{p247s1}), since
$Q(1-k+L)-1+\sum_{j=1}^{i}m_{\l_j} + \sum_{j=i+1}^{Q}(-m_{\bar\l_j}-1)
=-1$, we have the condition:

$$\sum_{j=1}^{i}(-1-m_{\l_j}) +
\sum_{j=i+1}^{Q}(m_{\bar\l_j}) = Q(L-k)$$ for $i \leq i \leq Q$
\mbox{$\square$}

Now we have the two key lemmas.  We will use Lemma \ref{43} to
limit the repetitions of negative modes and Lemma \ref{mainlemma}
to limit the repetitions of nonnegative modes.

\begin{thmn}
\label{mainthm} Let $V$ be a $C_2$ co-finite vertex operator
algebra, and let $W$ be an weak $V$ module generated by $w \in W$,
and let $x^i \in X$ for $1 \leq i \leq k$. Then $W$ is spanned by
elements of the form
$$x^{1}_{-n_1} x^{2}_{-n_2} \cdots x^{k}_{-n_k} w$$ where $n_1
\geq n_2 \geq \cdots \geq n_k
> -L$.  In addition, if  $n_j>0$, then $n_j > n_{j'}$ for
$j<j'$, and if $n_j \leq 0$ then $n_j=n_{j'}$ for at most $Q-1$
indices, $j'$.
\end{thmn}

These conditions say the following: For an element of the spanning
set, all the modes are decreasing and strictly less than L. If the
modes are negative then they are strictly decreasing. If the modes
are nonnegative then they are not strictly decreasing. There may
be repeats of nonnegative modes, but there are at most $Q-1$
repetitions.  Here is a sort of picture of the mode restrictions:

$$\mbox{strictly decreasing}<0 \leq \mbox{Q-1 repetitions}< L$$

The statement of the theorem might be easier to understand in a
simpler form:

{\it Let $W$ be an weak $V$ module generated by $w \in W$, and let
$x^i \in X$ for $1 \leq i \leq k$.  Then $W$ is spanned by
elements of the form $$x^{1}_{-n_1} x^{2}_{-n_2} \cdots
x^{k}_{-n_k} w$$ where the modes appear in decreasing order, are
less than $L$, and only finitely many of the same modes may appear
in the expression of each spanning set.}

This statement doesn't provide all the details, but it does give a
clearer picture.

\pf There are three main parts to this proof.  In the first part,
we define a filtration for the module and determine the properties
of the filtration.  The second part is the meat of proof, in which
we use an induction argument and two of the lemmas to impose
repetition restrictions on our spanning set elements. In the third
part of the proof, we demonstrate our spanning set elements indeed
do span $M$, and they obey the stated repetition restrictions.

\noindent{\bf Part I: Filtration and Properties}

First we define a filtration on $W$ as follows:
$$W^{(0)} \subset W^{(1)} \subset \cdots \subset W^{(t)} \subset
\cdots \subset W$$

\noindent Here $W^{(t)} = \mbox{span} \{u^1_{-n_1} u^2_{-n_2}
\cdots u^s_{-n_s}w\}$ where $ \sum_{i=1}^{s} wt(u^i) \leq t$,
$u^i$ a homogenous element of $V$. Now $W=\bigcup_t W^{(t)}$ since
$W$ is generated by $w$.

Let $u= u^1_{-n_1} \cdots u^s_{-n_s}w \in W^{(t)}$. Consider what
happens where two modes are transposed:
\begin{eqnarray*}
\lefteqn{u^1_{-n_1} \cdots u^{i+1}_{-n_{i+1}} u^{i}_{-n_{i}}
\cdots u^s_{-n_s}w}\\
&=&u^1_{-n_1} \cdots u^s_{-n_s}w - u^1_{-n_1} \cdots
[u^{i}_{-n_{i}}, u^{i+1}_{-n_{i+1}}] \cdots u^s_{-n_s}w
\end{eqnarray*}

\noindent Now $u^1_{-n_1} \cdots [u^{i}_{-n_{i}},
u^{i+1}_{-n_{i+1}}] \cdots u^s_{-n_s}w \in W^{(t-1)}$, because by
Lemma \ref{41}

$$[ u^{i}_{-n_{i}} , u^{i+1}_{-n_{i+1}} ]=\sum\limits_{j \geq 0}
{{wt(u^i)} \atopwithdelims ()
j}(u^i_{j}u^{i+1})_{-n_{i}-n_{i+1}-j}$$

\noindent and $wt(u^i_{j}u^{i+1})= wt(u^i)+wt(u^{i+1})-j-1 <
wt(u^i) + wt(u^{i+1})$. So,

$$u^1_{-n_1} \cdots u^{i+1}_{-n_{i+1}} u^{i}_{-n_{i}} \cdots u^s_{-n_s}w =
u^1_{-n_1} \cdots u^s_{-n_s}w +v$$ where $v \in W^{(t-1)}$. This
allows us to reorder modes of a vector in $W^{(t)}$ without
changing the vector modulo $W^{(t-1)}$.

Using the same $u$ as above, where $u \in W^{(t)}$. Now consider
what happens when we replace $u^i$ by its representative modulo
$C_2(V)$.  Let $u^i=x^i+c$  where $x^i \in X$ and $c \in V \slash
C_2(V)$. Without loss of generality we can assume $c=v_{-2}y$.

\begin{eqnarray*}
\lefteqn{u^1_{-n_1} \cdots u^{i}_{-n_{i}} \cdots u^s_{-n_s}w} \\
&=& u^1_{-n_1} \cdots (x^i+c)_{-n_{i}} \cdots u^s_{-n_s}w \\
&=& u^1_{-n_1} \cdots x^{i}_{-n_{i}} \cdots u^s_{-n_s}w \\
&&+ u^1_{-n_1} \cdots c_{-n_{i}} \cdots u^s_{-n_s}w \\
&=& u^1_{-n_1} \cdots x^{i}_{-n_{i}} \cdots u^s_{-n_s}w \\
&&+ u^1_{-n_1} \cdots (v_{-2}y)_{-n_{i}} \cdots u^s_{-n_s}w \\
\end{eqnarray*}

Now $wt(v_{-2}y)=wt(v)+wt(y)+1$, but by lemma \ref{42}, we can
replace $(v_{-2}y)_{-n_{i}}$ by:
$$\sum\limits_{k \geq 0}(-1)^k
{{-2} \atopwithdelims () k}v_{-2-k}y_{-n_{i}+k} -\sum\limits_{i
\geq 0}(-1)^{k-2} {{-2} \atopwithdelims () k}y_{-2-n_{i}-k}v_{k}$$
where this only contributes $wt(y)+wt(v)$ to the filtration level,
which is less than the $wt(v)+wt(y)+1$ that $v_{-2}y$ contributes.
Thus we have
$$u^1_{-n_1} \cdots u^{i}_{-n_{i}} \cdots u^s_{-n_s}w =
u^1_{-n_1} \cdots x^j_{-n_{i}} \cdots u^s_{-n_s}w +b$$ where $b
\in W^{(t-1)}$.  This means that we can replace $u^i$ by its
representative modulo $C_2(V)$ without changing the original
vector $u$ modulo $W^{(t-1)}$.

The upshot of these observations is that we can reorder modes and
replace vectors modulo $C_2(V)$ with out increasing the filtration
of the vector.  In fact, under reordering and replacement the
vector stays the same modulo a lower filtration level.

\noindent{\bf Part II: Induction}

To proceed with the proof, we now use induction on the pairs
$(t,K)$ with the ordering $(t,K) >(t',K')$ if $t>t'$ or $t=t'$ and
$K>K'$. The inductive hypothesis is:

$$W^{(t)}=\mbox{span} \{x^{1}_{-n_1} x^{2}_{-n_2} \cdots
x^{n}_{-n_k} w: x^i \in X \mbox{ for } 1 \leq i \leq Q\}$$

\noindent where $n_1 \geq n_2 \geq \cdots \geq n_k >-L$ and,
$\sum_i wt(x^i) \leq t$.  In addition,  if $K \geq n_j>0$ then
$n_j>n_{j+1}$. Also, for $n_j \leq 0$ and $n_j <K$ then $n_j =
n_{j'}$ for at most $Q-1$ indices, $j'$. Basically the induction
hypothesis is saying that there are repetition restrictions on
modes greater than $-K$, but no restrictions on modes less than
$-K$. The restrictions on modes greater than $-K$ are the
following: Negative modes greater than $-K$ are strictly
increasing, and positive modes greater than $-K$ can have at most
$Q-1$ repetitions.  This picture gives an idea of what is going on
with the modes in the induction hypothesis:

$$\mbox{decreasing} \leq -K < \mbox{strictly decreasing}  <0 \leq
\mbox{at most Q-1 repetitions }  < L$$

We proceed in the following manner:  We start with the base case
$(0,-L)$, then show the induction hypothesis holds for all pairs
$(t,-L)$, and then finally move on to the general case, $(t,K)$.

Our base case is $(0,-L)$, that is there are no repetition
restrictions for any modes. Let $u=u^1_{-n_1} \cdots u^s_{-n_s}w
\in W^{(0)}$. Since $V$ is CFT, $u^i=c_i {\textbf 1}$ where $c_i \in
\C$. If $-n_i \ne -1$ then $u=0$, because ${\textbf 1}_p=Id_V$ if and
only if $p=-1$.  If $u \ne 0$ then $u=c({\textbf 1}_{-1})^s w=cw$ , $c
\in \C$.  And $w$ is in our set of spanning elements.  So this
case is finished.

Now consider the pairs $(t,-L)$.  This our intermediate step in
the induction proof.  This step will show that we can reorder and
replace modes so that we can rewrite any vector as the sum of
elements of the form $x^1_{-n_1} x^2_{-n_2} \cdots x^l_{-n_l}w$
where the modes are decreasing.  Let $u=u^1_{-n_1} \cdots
u^k_{-n_k}w \in W^{(t)}$.  By using the results from Lemmas
\ref{41} and \ref{42}, we can reorder the modes and replace $u^i$
by $x^{i}$, its representative in $X$, to obtain

$$u=x^{1}_{-1} \cdots x^{k}_{-n_k} w + v$$

\noindent where $n_1 \geq n_2 \geq \cdots \geq n_k$ and $v \in
W^{(t-1)}$. Now if $-n_s \geq L$ ,  then $x^{k}_{-n_k}w=0$.  This
means we are left with $v$ which is in $W^{(t-1)}$, and by the
induction hypothesis for $(t-1,-L)$, $v=\sum_{s \in
S}x^{s_1}_{-n_{s_1}} x^{s_2}_{-n_{s_2}} \cdots x^{s_q}_{-n_{s_q}}
w$ where $n_{s_1} \geq n_{s_2} \geq \cdots \geq n_{s_q} >-L$, and
$S$ is a finite index set.  So if $-n_s \geq L$, we have $u=v$,
and $v$ is a spanning set element for the induction hypothesis
$(t-1,-L)$.

If $-n_k<L$,

$$u=x^{1}_{-n_1} \cdots x^{k}_{-n_k} w + v$$

\noindent where $n_1 \geq n_2 \geq \cdots \geq n_k >-L$ and $v \in
W^{(t-1)}$. Again using the induction hypothesis for $(t-1,-L)$,
$v=\sum_{s \in S}x^{s_1}_{-n_{s_1}} x^{s_2}_{-n_{s_2}} \cdots
x^{s_q}_{-n_{s_q}} w$ where $n_{s_1} \geq n_{s_2} \geq \cdots \geq
n_{s_q} >-L$, and $S$ is a finite index set.  Then

$$u= x^{1}_{-n_1} \cdots x^{k}_{-n_k} w + \sum_{s \in
S}x^{s_1}_{-n_{s_1}} x^{s_2}_{-n_{s_2}} \cdots x^{s_q}_{-n_{s_q}}
w$$

\noindent where all the modes are decreasing.  This completes our
intermediate step to show that the induction hypothesis holds for
all pairs $(t,-L)$.

This intermediate step ensures our spanning set has decreasing
modes.  The next step is to show that the modes greater than $-K$
can only repeat an finite number of times.  It is now that we move
to the general case, where we try to show that the induction
hypothesis hold for all pairs $(t,K)$.

We begin by assuming that the inductive hypothesis holds for all
pairs strictly less than $(t,K)$. Let $u=u^1_{-n_1} \cdots
u^k_{-n_k}w \in W^{(t)}$.  There are two cases; $K \leq 0$ and
$K>0$.  For the case where $K>0$, this corresponds to placing
repetition restrictions on negative modes, and we use Lemma
\ref{43} to do this.  When $K \leq 0$, we use Lemma
\ref{mainlemma} to impose repetition restriction on nonnegative
modes.

\noindent{\bf Case 1: $K \leq 0$}

In this case we are working on modes between $0$ and $L$, trying
to impose repetition restrictions. Again, let $u=u^1_{-n_1} \cdots
u^k_{-n_k}w \in W^{(t)}$. By rearranging and replacing terms by
Lemmas \ref{41} and \ref{42} ,and by applying the inductive
hypothesis for $(t,K-1)$, we get $u$ is the sum of vectors of the
form:
\begin{eqnarray}
\label{glob}\dot{x}^{1}_{-n_1} \cdots \dot{x}^{m}_{-n_m} x^1_{-K} \cdots x^p_{-K}
y^{1}_{-l_1} \cdots y^{s}_{-l_s}w+v
\end{eqnarray}
where $\dot{x}, x, y \in X$ for all indices; $n_1 \geq \cdots \geq
n_m > K > l_1 \geq \cdots \geq l_s >-L$;  $l_i=l_{i'}$ for at most
$Q-1$ indices; and $v \in W^{t-1}$. If $p<Q$, that is we have less
that $Q$ modes that are $-K$, then

$$\dot{x}^{1}_{-n_1} \cdots \dot{x}^{m}_{-n_m} x^1_{-K} \cdots
x^p_{-K} y^{1}_{-l_1} \cdots y^{s}_{-l_s}w$$

\noindent is in the proper form.  That is, it is a spanning set
element for $W^{(t)}$, as given in the induction hypothesis. We
finish this subcase off by applying the induction hypothesis for
$(t-1,K)$ to $v$ to obtain:

$$v=\sum_{s \in S}x^{s_1}_{-n_{s_1}} x^{s_2}_{-n_{s_2}} \cdots
x^{s_q}_{-n_{s_q}} w$$

\noindent where $n_{s_1} \geq n_{s_2} \geq \cdots \geq n_{s_q}
>-L$; $n_{s_j}=n_{s_{j+1}}$ is only allowed for $n_{s_j}>K$ and
$n_{s_j} > 0$; for $n_{s_j} \leq 0$ and $n_{s_j} <-K$ then
$n_{s_j} = n_{s_{j'}}$ for at most $Q-1$ indices, $j'$; and $S$ is
a finite index set.  Thus

$$u=\dot{x}^{1}_{-n_1} \cdots \dot{x}^{m}_{-n_m} x^1_{-K} \cdots
x^p_{-K} y^{1}_{-l_1} \cdots y^{s}_{-l_s}w+\sum_{s \in
S}x^{s_1}_{-n_{s_1}} x^{s_2}_{-n_{s_2}} \cdots x^{s_q}_{-n_{s_q}}
w$$

\noindent and $u$ is the sum of spanning set elements with mode
restrictions fitting the induction hypothesis for $(t,K)$.

Note that for the rest of the cases we deal with in this proof, we
can place $v$ in the desired form using the same method as above.
That is we can apply the induction hypothesis for $(t-1,K)$ to $v$
to place it in the desired form.  This is when $v$ is a vector
with filtration level $t-1$ that comes from reordering and
replacing the modes of a vector with filtration level $t$. Because
of this fact, in the remaining part of the proof we will not worry
about this $v$ that comes from reordering and replacing.

If $p \geq Q$ and $m\neq 0$, We have more than $Q$ modes that are
$-K$ in our expression.  We need to find a way to reduce the
number of repetitions of the mode $-K$ to $Q-1$. Since
$\sum_{i=1}^m wt(\dot{x}^i)>0$, we can apply the induction
hypothesis for $(t-1,K)$ to $x^1_{-K} \cdots x^p_{-K} y^{1}_{-l_1}
\cdots y^{s}_{-l_s}w$. When we do this we obtain a sum of vectors
of the form

$$\dot{x}^{1}_{-n_1} \cdots \dot{x}^{m}_{-n_m} {x'}^1_{-K} \cdots {x'}^{p'}_{-K}
y^{1}_{-l_1} \cdots y^{s}_{-l_s}w$$

\noindent where $p'<Q,$ and thus satisfies the proper repetition
restrictions.  The key here is that there are modes less than $-K$
at the front of this vector. Since each $x \in X$ has a positive
weight.  We can use the induction case for $(t-1,K)$ to the back
part of this vector to give the modes $-K$ and higher the proper
repetition restrictions.  We will use this method for reduction a
few more times.

If $m=0$, the we are dealing with a vector of the form

$$x^{1}_{-K} \cdots x^{p}_{-K}y^{1}_{-l_1} \cdots y^{s}_{-l_s}w$$

\noindent with $p\geq Q$.  Now we apply Lemma \ref{mainlemma}
where $L-k=-K$ to get:

\begin{eqnarray}
\lefteqn{(\prod_{j=1}^Q
x_{-K}^{j}) x^{Q+1}_{-K} \cdots x^p_{-K} y^{1}_{-l_1} \cdots y^{s}_{-l_s}w=}\\
&& (Res_z \{Y(\sum_{r \in
R} x^{r_1}_{-n_{r_1}} x^{r_2}_{-n_{r_2}} \cdots
x^{r_l}_{-n_{r_l}}{\textbf 1},z)\cdot z^{Q(K-1)}\}) \\
&& \cdot x^{Q+1}_{-K} \cdots x^p_{-K} y^{1}_{-l_1} \cdots y^{s}_{-l_s}w \\
&-& (\sum_{i=1}^{Q} \sum_{\l \in \L^i_Q} \sum_{m_{i} \geq 0} (\prod_{j=1}^{i}
x^{\l_j}_{-1-m_{\l_j}}) \label{sum1}(\prod_{j=i+1}^{Q} x^{\bar{\l}_j}_{m_{\bar\l_j}})) \\
&&\cdot x^{Q+1}_{-K} \cdots x^p_{-K} y^{1}_{-l_1} \cdots y^{s}_{-l_s}w \\
&-& (\sum_{m_{j} \geq 0, 1 \leq j \leq Q} x^{Q}_{m_{Q}}\cdots x^{1}_{m_{1}})
x^{Q+1}_{-K} \cdots x^p_{-K}  y^{1}_{-l_1} \cdots y^{s}_{-l_s}w \label{sum2}
\end{eqnarray}

Where $\sum_{j=1}^{i}(-m_{\l_j}-1) + \sum_{j=i+1}^{Q}m_{\bar\l_j}
= -QK$, in (\ref{sum1}). And $m_{j} \neq -K$ for some $j$ and
$\sum_{j=1}^{Q}m_{j}=-QK$, in (\ref{sum2}).

Let's look at the following term.$$-(\sum_{m_{j} \geq 0, 1 \leq j
\leq Q} x^{Q}_{m_{Q}}\cdots x^{1}_{m_{1}}) x^{Q+1}_{-K} \cdots
x^p_{-K}  y^{1}_{-l_1} \cdots y^{s}_{-l_s}w$$ In this term,
$\sum_{j=1}^{Q}m_{j}=-QK$ and $m_{j} \neq -K$ for some $j$.  We
need to show that in this term there is a mode less than -K.  If
there is we can reorder the modes so that there is a mode less
than $-K$ at the front of the vector.  As before we can finish off
this case by applying the induction hypothesis on $(t-1,K)$ to the
back of the vector where the modes are greater than or equal to
$-K$.  So ${1 \over Q} \sum_{j=1}^{Q}(m_{j})=-K$ where one of the
$m_{j} \neq -K$. If one of the modes is not $-K$, then one of the
modes is less than $-K$ or greater than $-K$.  If the mode is less
than $-K$, we are done.  If the mode is greater than $-K$, since
the average of the modes is $-K$ then there must be another mode
less than $-K$, and we can apply the induction hypothesis. So we
are done with this term.

Next we examine,

$$(\prod_{j=1}^{i} x^{\l_j}_{-1-m_{\l_j}}) (\prod_{j=i+1}^{Q}
x^{\bar{\l}_j}_{m_{\bar\l_j}})) x^{Q+1}_{-K} \cdots x^p_{-K}
y^{1}_{-l_1} \cdots y^{s}_{-l_s}w$$ In this term, since $i \neq
0$, there will a mode of the form $x^{\l_j}_{-1-m_{\l_j}}$. In
this mode, $-1-m_{\l_j}$ is negative, and thus it will be less
than $-K$.  Again by reordering we are in the case where there is
a mode less than $-K$ at the front of the vector.  So by the
method described previously we are done with this term.

The final term we need to examine for this case is:

$$(Res_z \{Y(\sum_{r \in R} x^{r_1}_{-n_{r_1}} x^{r_2}_{-n_{r_2}} \cdots
x^{r_l}_{-n_{r_l}}{\textbf 1},z) z^{Q(K-1)}\}) x^{Q+1}_{-K} \cdots
x^p_{-K}  y^{1}_{-l_1} \cdots y^{s}_{-l_s}w$$

\noindent In this term, we must determine what happens when we
evaluate

$$Res_z \{Y(\sum_{r \in R} x^{r_1}_{-n_{r_1}} x^{r_2}_{-n_{r_2}}
\cdots x^{r_l}_{-n_{r_l}}{\textbf 1},z) z^{Q(K-1)}\}$$

\noindent After evaluating the residue we will get a sum of
products of $r_l$ modes. We must examine what happens for each
product of $r_l$ modes.  If the product of modes has at least one
negative mode or at least one mode less than $-K$, we can place
this vector in the proper form by rearranging the modes and
applying the induction hypothesis for $(t-1,K)$ the back part of
this vector.

What now remains are the products of modes for which each mode is
greater than of equal to $-K$.  But since $l < Q$,

$$(Res_z \{Y(\sum_{r \in
R} x^{r_1}_{-n_{r_1}} x^{r_2}_{-n_{r_2}} \cdots
x^{r_l}_{-n_{r_l}}{\textbf 1},z) z^{Q(K-1)}\})x^{Q+1}_{-K} \cdots
x^p_{-K}  y^{1}_{-l_1} \cdots y^{s}_{-l_s}w$$

\noindent has strictly less modes equal to $-K$, then

$$(\prod_{j=1}^Q x_{-K}^{p+1-j}) x^{Q+1}_{-K} \cdots x^p_{-K}
 y^{1}_{-l_1} \cdots y^{s}_{-l_s}w$$

\noindent does.  So this means that we can repeat the process of
applying Lemma \ref{mainlemma}, now to terms in

$$(Res_z \{Y(\sum_{r \in R} x^{r_1}_{-n_{r_1}} x^{r_2}_{-n_{r_2}} \cdots
x^{r_l}_{-n_{r_l}}{\textbf 1},z) z^{Q(K-1)}\})x^{Q+1}_{-K} \cdots
x^p_{-K}  y^{1}_{-l_1} \cdots y^{s}_{-l_s}w$$

\noindent By repeating this process we eventually, reduce the
number of modes equal to $-K$ to $Q-1$.  This finishes off the
case when $K \leq 0$.

\noindent{\bf Case 2: $K>0$}

In this case we show that the negative modes must be strictly
decreasing.  Start with $u$ as above.  Using the inductive
hypothesis for $(t, K-1)$, we get $u$ is the sum of vectors of the
form

$$\dot{x}^{1}_{-n_1} \cdots \dot{x}^{m}_{-n_m} x^{1}_{-K} \cdots
 x^p_{-K}  y^{1}_{-l_1} \cdots y^{s}_{-l_s}w$$
where $n_1 \geq \cdots \geq n_m > K > l_1 \geq \cdots \geq l_s
>-L$; $l_i=l_{i+1}$ only if $l_i \leq 0$; and if $n_j \leq 0$ then
 $n_j = n_{j'}$ for at most $Q-1$ indices, $j'$.  Now if $m \ne 0$,
then we are in the case where the is a mode at the front of the
vector that is less than $-K$.  As we have done before, we can
apply the induction hypothesis for $(t-1,K)$ to
$$x^{1}_{-K} \cdots x^p_{-K}  y^{1}_{-l_1} \cdots y^{s}_{-l_s}w$$  This
will remove the repeats at $-K$, giving us a vector of the form
$$\dot{x}^{1}_{-n_1} \cdots \dot{x}^{m}_{-n_m}x^{1}_{-l_1} \cdots
 x^p_{-l_p}w$$ where $n_1 \geq \cdots \geq n_m > K \geq l_1 \geq \cdots \geq l_s
>-L$; $l_i=l_{i+1}$ only if $l_i \leq 0$; and if $n_j \leq 0$ then
$n_j = n_{j'}$ for at most $Q-1$ indices, $j'$.  These are the
mode restrictions for the induction hypothesis for $(t,K)$, so
this vector is in the required form.  Again this case where there
is a mode less than $-K$ at the front of the vector is an
important case, and will be used again.

So now we are reduced to considering what happens if $m= 0$.  That
is, we are dealing with a vector of the form
$$x^{1}_{-K} \cdots x^p_{-K}  y^{1}_{-l_1} \cdots y^{s}_{-l_s}w$$
where $K > l_1 \geq \cdots \geq l_s >-L$ and
$l_i=l_{i+1}$ only if $l_i \leq 0$.  Now we will use Lemma \ref{43}.

\begin{eqnarray}
\lefteqn{x^{1}_{-K} \cdots x^{p}_{-K} y^{1}_{-l_1}
\cdots y^{s}_{-l_s}w} \\
&=& (x^{1}_{-K} x^{2})_{-2K+1} x^{3}_{-K} \cdots
x^{p}_{-K}
y^{1}_{-l_1} \cdots y^{s}_{-l_s}w \\
&&- \sum\limits_{i > 0}(-1)^i
{{-K} \atopwithdelims () i}x^{1}_{-1-i}x^{2}_{-2K+1+i}
x^{3}_{-K} \cdots x^{p}_{-K}y^{1}_{-l_1}
\cdots y^{s}_{-l_s}w\\
&&+ \sum\limits_{i \geq 0} (-1)^{i-K}
{{-K} \atopwithdelims () i} x^{2}_{-2K-i}
x^{1}_{i} x^{3}_{-K} \cdots x^{p}_{-K}y^{1}_{-l_1}
\cdots y^{s}_{-l_s}w
\end{eqnarray}

In the first term, $2K-1>K$ if and only if  $K>1$.  We deal with
the case where $K=1$ in a moment.  For now assume $K>1$. In the
second term, it is easy to check that either $-1-i<-K$ or
$-2K+1+i<-K$ for all $i \geq 0$ and $i \neq K-1$. In the third
term, $-2K+1+i <-K$ for all $i \geq 0$. So now, we have reduced
this to the case where there is a mode at the front of each vector
that is less than $-K$.  We have shown previously that vectors of
this form can be place in the required form.

Now, we deal with the case where $K=1$.  This case requires us to
repeatedly apply Lemma \ref{43}.  We are dealing with a vector of
the form

$$x^{1}_{-1} \cdots x^{p}_{-1} y^{1}_{-l_1}
\cdots y^{s}_{-l_s}w$$

By applying Lemma \ref{43}, we obtain

\begin{eqnarray}
\lefteqn{x^{1}_{-1} \cdots x^{p}_{-1} y^{1}_{-l_1}
\cdots y^{s}_{-l_s}w} \\
&=& (x^{1}_{-1} x^{2})_{-1} x^{3}_{-1} \cdots
x^{p}_{-1}
y^{1}_{-l_1} \cdots y^{s}_{-l_s}w \label{neg1}\\
&&- \sum\limits_{i > 0}(-1)^i
{{-1} \atopwithdelims () i}x^{1}_{-1-i}x^{2}_{-1+i}
x^{3}_{-1} \cdots x^{p}_{-1}y^{1}_{-l_1}
\cdots y^{s}_{-l_s}w \label{neg2}\\
&&+ \sum\limits_{i \geq 0} (-1)^{i-1}
{{-1} \atopwithdelims () i} x^{2}_{-2-i}
x^{1}_{i} x^{3}_{-1} \cdots x^{p}_{-1}y^{1}_{-l_1}
\cdots y^{s}_{-l_s}w \label{neg3}
\end{eqnarray}

For (\ref{neg2}) and (\ref{neg3}), we have a mode at the front of
the vector which is strictly less than $-1$.  As before, we can
apply the induction hypothesis for $(t-1,1)$ the latter $s-1$
modes to place these vectors in the required form. For
(\ref{neg1}), we can replace $(x^{k_1}_{-1} x^{k_2})_{-1}$ by its
representative modulo $C_2(V)$.  The we have an expression with on
less $-1$ mode.  By repeating the process, we eventually eliminate
all repetitions of the $-1$ mode.

\noindent{\bf Part III: A Spanning Set}

Finally, we must show that given any element of $W$ we can
rewritten it a sum of spanning set elements, as claimed.  We know
that $M$ is spanned by elements of the form

$$u^1_{-m_1} u^2_{-m_2} \cdots u^k_{-m_k} w$$

\noindent We must show that this vector can be written as a sum of
elements of the form

$$x^{1}_{-n_1} x^{2}_{-n_2} \cdots
x^{l}_{-n_l} w$$

\noindent where $n_1 \geq n_2 \geq \cdots \geq n_k
> -L$.  In addition, if  $n_j>0$, then $n_j > n_{j'}$ for
$j<j'$, and if $n_j \leq 0$ then $n_j=n_{j'}$ for at most $Q-1$
indices, $j'$.

We can make the assumption that $u^i$ is homogenous for $1 \leq i
\leq k$.  Any vector is a sum of homogenous vectors and
$(u+v)_n=u_n+v_n$ where $u,v \in V$ and $n \in \Z$.  Let
$D=wt(u^1_{-m_1} u^2_{-m_2} \cdots u^k_{-m_k})=\sum_{i=1}^k
wt(u^i)+m_i-1$ and let $t=\sum_{i=1}^k wt(u^i)$, the filtration
level of $u^1_{-m_1} u^2_{-m_2} \cdots u^k_{-m_k} w$.  Consider
the induction hypothesis for the pair $(t,D+(t-1)L+1)$.  By the
induction hypothesis for this pair we have

$$u^1_{-m_1} u^2_{-m_2} \cdots u^k_{-m_k} w= \sum_{r \in R}
x^{r_1}_{-n_{r_1}} \cdots x^{r_l}_{-n_{r_l}}w$$

\noindent where $n_{r_1} \geq n_{r_2} \geq \cdots \geq n_{r_l}
>-L$ , $\sum_i wt(x^i) \leq t$ , and $n_{r_j}=n_{r_{j+1}}$ is only
allowed for $n_{r_j}>D+(t-1)L+1$ and $n_j > 0$.  Also, for $n_j
\leq 0$ and $n_j <D+(t-1)L+1$ then $n_{r_j} = n_{r_{j'}}$ for at
most $Q-1$ indices, $r_{j'}$.  Note that
$$wt(u^1_{-m_1} u^2_{-m_2} \cdots u^k_{-m_k})= wt(x^{i_1}_{-n_{r_1}}
\cdots x^{i_l}_{-n_{r_l}})$$  This is true because the identities
that we use to prove the induction hypothesis all preserve the
weights of operators.   Now we have that

\begin{eqnarray}
D &=& wt(x^{r_1}_{-n_{r_1}} \cdots x^{r_l}_{-n_{r_l}})\\
&=&\sum^l_{i=1}(wt(x^{r_i})+n_{r_i}-1)
\end{eqnarray}

\noindent Now we attempt to calculate an lower bound for
$n_{r_1}$.

\begin{eqnarray}
n_{r_1} &=& D -wt(x^{r_1})-1-\sum^l_{i=2}(wt(x^{r_i})+n_{r_i}-1)\\
& \leq& D-\sum^l_{i=2}n_{r_i}\\
&\leq& D+(l-1)L
\end{eqnarray}
But $\sum_{i=1}^l wt(x^{r_i}) \leq t$ implies that $l \leq t$,
since $wt(x^i) >0$.  So we have
$$n_{r_1} \leq D+(t-1)L$$

\noindent Since we used the induction hypothesis for the pair
$(t,D+(t-1)L+1)$, we have

 $$u^1_{-m_1} u^1_{-m_2} \cdots
u^1_{-m_k} w= \sum_{r \in R} x^{r_1}_{-n_{r_1}} \cdots
x^{r_l}_{-n_{r_l}}w$$

\noindent where $n_{r_1} \geq n_{r_2} \geq \cdots \geq n_{r_l}
>-L$. In addition, if $n_j>0$, then $n_{r_j}>n_{r_{j'}}$ for
$j<j'$ and if $n_j \leq 0$ then $n_j=n_{j'}$ for at most $Q-1$
indices, $j'$. This show that $u^1_{-m_1} u^1_{-m_2} \cdots
u^1_{-m_k} w$ can be rewritten in terms of spanning set vectors of
the desired form. \mbox{$\square$}

So this theorem accomplishes the goal of imposing a finite repeat
condition on the modes of the module spanning set.  This tells us
that limiting the number of positive modes, and thus we limit the
number of modes with negative weights.  The idea here is that
modes with negative weights push our vector $w$ down.  Limiting
the number of positive modes means we can only push $w$ down so
far.  Now since there is no grading on weak modules, this is not a
precise statement, but it is the picture of what is going on.

\section{Additional Results}

In this section we look at some applications of this module
spanning set.  So all of these results assume the $V$ is $C_2$
co-finite.

In the work of Li \cite{L}, he defines $C_n(M)$ and used it to
show that the fusion rules are finite for admissible modules.

\begin{de}
Let $M$ be a weak $V$ module, then $C_n(M)=\{v_{-n}m \mid m\in M, v \in
V\}$
\end{de}

As with vertex operator algebras, we can look at the quotient
space $M \slash C_n(M)$.

\begin{de}
A weak $V$ module, $M$, is called $C_n$ co-finite if $M \slash
C_n(M)$ is finite dimensional.
\end{de}

Just as $C_2$ co-finiteness implies $C_n$ co-finiteness for $V$,
the $C_2$ co-finiteness implies the co-finiteness of $C_n(M)$.

\begin{cor}
If $M$, a irreducible weak  $V$ module, is $C_2$ co-finite for
$n>2$ then $M$ is $C_n$ co-finite for $n \geq 2$.
\end{cor}

\pf Using our new modules spanning set, we see that $M \slash
C_n(M)$ is spanned by elements of the form $$x^{1}_{-n_1} x^{2}_{-n_2} \cdots
x^{k}_{-n_k} w$$ where $n >n_1 \geq n_2 \geq \cdots \geq n_k
> -L$.  In addition, if  $n_j>0$, then $n_j > n_{j'}$ for
$j<j'$, and if $n_j \leq 0$ then $n_j=n_{j'}$ for at most $Q-1$
indices, $j'$.  Since $M$ is $C_2$ co-finite, and there are only
finitely many vectors that we add to a spanning set of $M \slash
C_2(M)$ to get a spanning set for $M \slash C_n(M)$, then $M$ is
$C_n$ co-finite. \mbox{$\square$}

In \cite{DLM2}, they show that if $V$ is $C_2$ co-finite then
$A(V)$ is finite dimensional.  Using the module spanning set, we
can extend this result to cover $A_n(M)$.  The $A(V)$ bimodule
$A(W)$ appears first in \cite{FZ}.  We can extend the definition
of $A(W)$ to $A_n(W)$, just as the definition of $A(V)$ is
extended to $A_n(V)$ in \cite{DLM2}.

\begin{de}
Let $u \circ_n w = Res_z Y(u,z)w{{(1+z)^{wt(u)+n}}
\over{z^{2n+2}}}$.  For homogenous $u$, this can be rewritten as
$u \circ_n w =\sum_{j \geq 0} {{wt(u)+n} \atopwithdelims () j}
u_{j-2n-2}w$ where $u \in V$ and $w \in M$. Then let $O_n(M)= \{ u
\circ_n w : u \in V, w \in M \}$.  Finally we define $A_n(M)= M
\slash O_n(M)$.  We use the notation convention that $A_0(M)=A(M)$
\end{de}

\begin{cor}
\label{c55} If $M$,a irreducible weak $V$ module, is $C_2$
co-finite then $A_n(M)$ is finite dimensional.
\end{cor}

\pf First if $V$ is $C_2$ co-finite then $M$ is $C_{2n+2}$ finite.
Let $u$ be a module spanning set element, that is $u=x^1_{-n_1}
x^2_{-n_2} \cdots x^k_{-n_k}w$.  Now define $N$ to be the smallest
integer such that if $wt(x^1_{-n_1} x^2_{-n_2} \cdots x^k_{-n_k})
\geq N$ then $n_1 \geq 2n+2$.  We will show that each element of $M$ can be
written as an element in $O_n(M)$ plus an element in a finite
dimensional space. Specifically we will show for any module
spanning element $u$ that $u=x+y$ where $y \in O_n(M)$ and
$x=\sum_{r \in R} x^{r_1}_{-n_{r_1}} \cdots x^{r_l}_{-n_{r_l}}w$
where $wt(x^{r_1}_{-n_{r_1}} \cdots x^{r_l}_{-n_{r_l}}) < N$. That
is, $x$ is in a finite dimensional subspace of $M$

We prove this by induction on $r=wt(x^1_{-n_1} x^2_{-n_2} \cdots
x^k_{-n_k})-N$.  If $r \leq 0$ then $u=x+y$ where $y=0$.  Now
consider $u=x^1_{-n_1} x^2_{-n_2} \cdots x^k_{-n_k}w$ where
$wt(x^1_{-n_1} x^2_{-n_2} \cdots x^k_{-n_k})-N=r>0$. This means
that $u \in C_{2n+2}(M)$.  So, $n_1 \geq 2n+2$.  We use the fact
that $(L(-1)v)_{-n}=n v_{-n-1}$, we can rewrite $x^1_{-n_1}$ as
$({1 \over s!} L(-1)^s x^1)_{-2n-2}$ for some positive $s$.  So
now

\begin{eqnarray}
u &=& x^1_{-n_1} x^2_{-n_2} \cdots x^k_{-n_k}w \\
&=& ({1 \over s!} L(-1)^s x^1)_{-2n-2}x^2_{-n_2} \cdots
x^k_{-n_k}w \\
&=& ({1 \over s!} L(-1)^s x^1) \circ_n  (x^2_{-n_2}
\cdots x^k_{-n_k}w) \\
&&- \sum_{j \geq 1} {{wt(L(-1)^s x^1)+n} \atopwithdelims () j} ({1
\over s!} L(-1)^s x^1)_{j-2n-2} x^2_{-n_2} \cdots x^k_{-n_k}w\\
&=& ({1 \over s!} L(-1)^s x^1) \circ_n  (x^2_{-n_2}
\cdots x^k_{-n_k}w) \\
&&-\sum_{r \in R} x^{r_1}_{-n_{r_1}} \cdots x^{r_l}_{-n_{r_l}}w
\end{eqnarray}
where $wt(x^{r_1}_{-n_{r_1}} \cdots x^{r_l}_{-n_{r_l}})-N< r$. In
the last step we rewrite $$\sum_{j \geq 1} {{wt(L(-1)^s x^1)+n}
\atopwithdelims () j} ({1 \over s!} L(-1)^s x^1)_{j-2n-2}
x^2_{-n_2} \cdots x^k_{-n_k}w$$ in terms of the module spanning
set. We know that $({1 \over s!} L(-1)^s x^1) \circ_n  (x^2_{-n_2}
\cdots x^k_{-n_k}w) \in O_n(M)$. We can apply induction argument
to $$\sum_{r \in R} x^{r_1}_{-n_{r_1}} \cdots x^{r_l}_{-n_{r_l}}
w$$ so that we can place it in the form $x+y$ as described above.
Now since $A_n(M)=M\slash O_n(M)$, $A_n(M)$ is spanned by these
elements $x+O_n(M)$ where $x$ in in a finite dimensional subspace
of $M$.\mbox{$\square$}

\begin{cor}
Let  $W=\bigoplus_{n \geq 0}W(n)$ be an admissible $V$ module generated by a finite dimensional $W(i)$, then $W$ is an ordinary $V$ module.
\end{cor}

\pf To prove this Corollary, it is sufficient to show that each
graded piece of $W$ is finite dimensional.  Let $\{w^j\}$ with $1
\leq j \leq s$ be a basis for $W(i)$. For each $w^j$ we can look a
the module spanning set associated to $w^j$. Since $W(i)$ generates $W$, these
spanning sets combined will form a spanning set for $W$. Since $W$ is
admissible, if $wt(x^1_{-n_1}x^2_{-n_2} \cdots x^k_{-n_k})=N$ then
$x^1_{-n_1}x^2_{-n_2} \cdots x^k_{-n_k}w^j \in W(i+N)$. This means
that $W(n)$ will spanned by modules spanning set elements of the
form:

$$x^1_{-n_1}x^2_{-n_2} \cdots x^k_{-n_k}w^j$$

\noindent where $wt(x^1_{-n_1}x^2_{-n_2} \cdots x^k_{-n_k})=n-i$.
Because of the mode restrictions, this spanning set will be
finite. \mbox{$\square$}


\begin{thebibliography}{ABCD}


\bibitem[B]{B}
R. Borcherds., \emph{Vertex algebras, Kac-Moody algebras, and the
Monster}, Proc. Natl. Acad. Sci. USA {\bf 83}, 1986, 3068-3071.

\bibitem[D]{D}
C.Dong., \emph{Introduction to vertex operator algebras I.},
arXiv:q-alg/9504017.

\bibitem[DLM]{DLM}
C.Dong, H. Li and G. Mason., \emph{Regularity of rational vertex
operator algebra}, Adv. Math. {\bf 312}, 1997, 148-166.

\bibitem[DLM2]{DLM3}
C.Dong, H. Li and G. Mason., \emph{Twisted representations of
vertex operator algebras}, Math. Ann. {\bf 310}, 1998, 571-600.

\bibitem[DLM3]{DLM2}
C.Dong, H. Li and G. Mason., \emph{Vertex operator algebras and
associative algebras}, J. Algebra {\bf{206}}, 1998, 67-96.

\bibitem[FF]{FF}
B. Feigin, and D. Fuchs., \emph{Verma modudules over the Virasoro
algebra}, Lecture Notes in Math. {\bf{1060}}, 1984, 230-245.

\bibitem[FHL]{FHL}
I. Frenkel, Y.-Z. Huang, and J. Lepowsky., \emph{On axiomatic
approaches to vertex operator algebras and modules}, Mem. Amer.
Math. Soc.{\bf{104}}, 1993.


\bibitem[FLM]{FLM}
I. B. Frenkel, J. Lepowsky, and A. Meurman., \emph{ Vertex
Operator Algebras and the Monster}, Pure and Appl. Math., {\bf
Vol. 134}, Academic Press, Boston 1988.

\bibitem[FZ]{FZ}
I. Frenkel, and Y.-C. Zhu., \emph{Vertex operator algebras
associated to representations of affine and virasoro algebras},
Duke Math J. {\bf{66}}, 1992, 123-168.

\bibitem[GN]{GN}
M. Gaberdiel, and A. Neitzke, \emph{Rationality, quasirationality,
and finite W-algebras}, DAMTP-200-111.


\bibitem[KL]{KL}
H. Li, and M. Karel., \emph{Certain generating subspace for vertex
operator algebras}, Adv. Math. {\bf 132}, 1997, 148-166.


\bibitem[L]{L}
H. Li., \emph{Some finiteness properties of regular vertex
operator algebras}, J. Algebra. {\bf 212}, 1999, 495-514.


\bibitem[L2]{L2}
H. Li., \emph{Determining fusion rules by a $A(V)$-modules and
bimodules.}, J. Algebra. {\bf 212}, 1999, 515-556.


\bibitem[Z]{Z}
Y. Zhu., \emph{ Modular invariance of characters of vertex
operater algebras}, J. Amer. Math. Soc. {\bf 9}, 1996, 237-302.

\end{thebibliography}
\end{document}